\documentclass{amsart}

\usepackage{amsmath}
\usepackage{amsthm}
\usepackage{amssymb}
\usepackage{amsfonts}
\usepackage{bm}
\usepackage[shortlabels]{enumitem}
\usepackage[bookmarks=true,hyperindex,pdftex,colorlinks,citecolor=red, linkcolor=cyan]{hyperref}
\usepackage{centernot} 
\usepackage{marginnote} 
\usepackage{bbm}
\usepackage[all]{xy}
\usepackage{tikz}
\usetikzlibrary{arrows}
\usetikzlibrary{shapes,decorations}
\usetikzlibrary{positioning}
\usepackage{tikz-cd}
\usepackage[normalem]{ulem}
\usepackage{array}
\usepackage{bbold}

\allowdisplaybreaks

\newcommand{\C}{\mathbb{C}}             

\newcommand{\norm}[1]{\left\|{#1}\right\|}                  

\renewcommand{\restriction}{\mathord{\upharpoonright}}

\def\<{\langle}
\def\>{\rangle}

\renewcommand{\leq}{\leqslant}
\renewcommand{\geq}{\geqslant}

\newcommand{\fonction}[5]{\begin{array}{l rcl}
#1: & #2 & \longrightarrow & #3 \\
      & #4 & \longmapsto & #5 \end{array}}

\theoremstyle{plain}
\newtheorem{theorem}{Theorem}[section]
\newtheorem{lemma}[theorem]{Lemma}
\newtheorem{corollary}[theorem]{Corollary}
\newtheorem{proposition}[theorem]{Proposition}

\theoremstyle{definition}
\newtheorem*{definition*}{Definition}
\newtheorem{definition}[theorem]{Definition}
\newtheorem{example}[theorem]{Example}
\newtheorem{question}{Question}

\theoremstyle{remark}
\newtheorem{remark}[theorem]{Remark}

\begin{document}

\title{On universal operators on Banach spaces}

\author[M. Lemay]{Mathis Lemay}
\address[M. Lemay]{Univ Gustave Eiffel, Univ Paris Est Creteil, CNRS, LAMA UMR8050, F-77447 Marne-la-Vallée, France}
\email{mathis.lemay@univ-eiffel.fr}

\subjclass[2020]{47A15, 47B37, 46B20}
\keywords{Universal operators, invariant subspaces, Caradus criterion, spectrum of an operator}

\begin{abstract}
The paper studies universal operators and $n$-tuples of operators on Banach spaces. 
In the Hilbert-space setting, sufficient criteria for universality, concrete examples, and spectral properties are well established. We investigate extensions of these results to Banach spaces. In particular, we discuss the link between universality and the invariant subspace problem: this leads us to define a stronger notion of universality which is the notion of hereditarily universal operators. We give sufficient conditions for an operator to be hereditarily universal and provide some examples. We also study $n$-tuples of operators by generalizing Vladimir Müller's results to the Banach setting. Finally, we obtain a result on the essential spectrum of a universal $n$-tuple of operators on $\ell_p$.
\end{abstract}

\maketitle

\section{Introduction and preliminaries}

Let $X$ be a complex Banach space. By an \emph{operator on $X$}, we mean a bounded linear map $T : X \longrightarrow X$ and we write $\mathcal{L}(X)$ for the algebra of bounded operators from $X$ to $X$. For complex Banach spaces $X$ and $Y$, we say that $T \in \mathcal{L}(X)$ and $S \in \mathcal{L}(Y)$ are \emph{similar} if there exists a bounded linear isomorphism $J : X \longrightarrow Y$ with bounded inverse such that $SJ=JT$, or equivalently $T=J^{-1}SJ$. The notation $X \simeq Y$ means that $X$ and $Y$ are linearly isomorphic. The main notion addressed in this paper is the notion of \emph{universal operators} which is defined as follows:

\begin{definition}[Rota \cite{Rota1959}]\label{defuniv}
Let $X$ be an infinite-dimensional Banach space and let $U\in\mathcal L(X)$.
We say that $U$ is \emph{universal on $X$} (or \emph{universal in the sense of Rota})
if for every $T\in\mathcal L(X)$ there exist a closed $U$-invariant subspace
$M\subset X$ (that is, $U(M)\subset M$) and a scalar $c\in\C\setminus\{0\}$ such that
the operators $U\restriction_M:M\to M$ and $cT:X\to X$ are similar.
\end{definition}

In other words, a universal operator ``encodes" every operator of the space up to a multiplicative constant and a restriction to an invariant subspace. The definition is often presented first for a separable infinite-dimensional Hilbert space $H$. In that setting, the particular model of $H$ is irrelevant: any two such Hilbert spaces are isometric, so universality is an intrinsic property of the similarity class of the operator $U$.

The interest of studying such operators on Hilbert spaces lies in the link with the famous \emph{invariant subspace problem}. General background on universal operators and their role in invariant subspace theory can be found in \cite[Chapter 8]{CP2011} and in the survey \cite{CGGintro}, but we discuss this connection in more detail in Subsection~\ref{invsub}.

In 1969, Caradus provided elegant sufficient conditions for an operator to be universal on an infinite-dimensional separable Hilbert space. In Subsection \ref{Caradusetex}, we recall the statement of this famous theorem, and some well-known examples of universal operators on Hilbert spaces. The wish to establish a link between universal operators and the invariant subspace problem on Banach spaces leads us to define a stronger notion: hereditarily universal operators, that is, operators that model, up to similarity and a non-zero scalar factor, every operator on every infinite-dimensional closed subspaces. We show a similar ``Caradus-type" theorem to provide sufficient conditions for an operator to be hereditarily universal. Then, some geometrical properties of infinite-dimensional Banach spaces, like minimality or primeness, will allow us to make these sufficient conditions even simpler. This will lead to examples on sequence spaces. Then, we will exhibit spaces with no universal operators. 

The next part of this paper is dedicated to the study of the spectrum of universal operators on $\ell_p$. By generalizing the proof provided by Rikka Schröderus and Hans-Olav Tylli in \cite{schro}, we show that the spectrum contains a whole ball centered at $0$, and that every point of this ball is an eigenvalue with infinite multiplicity. Moreover, for every point of this ball, there exists an eigenvector which depends holomorphically on the eigenvalue.

Then we focus on $n$-tuples of operators. Some conditions on the kernel of the operators, combined with the existence of a bounded below operator, ensure that an $n$-tuple is universal. Finally, we give an idea of the elements of the essential spectrum of a universal $n$-tuple of operators on $\ell_p$.

\subsection{The relevance of universal operators regarding the invariant subspace problem}\label{invsub}

\begin{definition}
Let $X$ be a Banach space and $T\in\mathcal L(X)$. A closed subspace $M\subset X$ is called
\emph{$T$-invariant} if $T(M)\subset M$. It is called \emph{non-trivial} if $M\neq\{0\}$ and
$M\neq X$.

A non-trivial closed $T$-invariant subspace $M$ is \emph{minimal} if the only closed $T$-invariant subspaces contained in $M$ are $\{0\}$ and $M$.
\end{definition}

The classical \emph{invariant subspace problem} asks whether every bounded operator on a separable
infinite-dimensional Hilbert space has a non-trivial closed invariant subspace; this problem is open.
For general Banach spaces the answer is negative: there exist operators with no non-trivial invariant
subspaces (see, e.g., Enflo's early announcement \cite{Enflo} and Read's constructions
\cite{Read}, \cite{Read2}). However, the problem is still open for reflexive Banach spaces. For background we refer to \cite{CP2011}.

In the Hilbert setting, the link between universal operators and invariant subspaces can be summarized as follows.

\begin{proposition}\label{prop:ISP}
Let $H$ be an infinite-dimensional separable Hilbert space and let $U\in\mathcal L(H)$ be universal on $H$.
Then the following properties are equivalent:
\begin{enumerate}[(i)]
\item every operator $T\in\mathcal L(H)$ has a non-trivial invariant subspace;
\item every minimal non-trivial $U$-invariant subspace has finite dimension;
\item every minimal non-trivial $U$-invariant subspace has dimension $1$.
\end{enumerate}
\end{proposition}

\begin{proof}
We first note that $(iii) \implies (ii)$ is obvious.

$(ii) \implies (iii)$: 
Let $M\subset H$ be a minimal non-trivial $U$-invariant subspace and assume
$\dim M<\infty$. Since we work over $\C$, the operator $U\restriction_M$ has an eigenvalue and hence
a $1$-dimensional invariant subspace. This eigenline is a nonzero $U$-invariant subspace contained in $M$, so minimality forces it to equal $M$.

$(i) \implies (iii)$:
Let $M$ be a minimal non-trivial $U$-invariant subspace. If $\dim(M)$ is finite, then we are done because $U\restriction_{M}$ has an eigenvector, whose span is a one-dimensional invariant subspace. Minimality therefore forces $\dim(M)=1$. If $\dim(M)$ is infinite, then let $J : H \longrightarrow M$ be an isometric isomorphism, and let $T=J^{-1}(U\restriction_{M})J \in \mathcal{L}(H)$. By $(i)$, $T$ has a non-trivial closed invariant subspace $N \subset H$. Then $J(N)$ is a non-trivial closed $U$-invariant subspace of $M$, contradicting minimality. 
\medskip

$(iii) \implies (i)$: 
Fix $T\in\mathcal L(H)$ and suppose that $T$ has no non-trivial invariant subspace. Since $U$ is universal, there exist a closed
$U$-invariant subspace $M\subset H$, a scalar $c\neq 0$ and an isomorphism $J:H\longrightarrow M$ such that
$$
U\restriction_M \, J = J \, (cT).
$$
The subspace $M$ is proper in $H$: otherwise $U$ would be similar to $cT$, whereas universality applied to the zero operator forces $\ker(U)$ to contain a copy of $H$, and this kernel is a nonzero proper invariant subspace whenever $U\neq0$. Also $U\neq0$, since it models $I_H$. Every nonzero closed invariant subspace of $M$ pulls back to a nonzero closed $T$-invariant subspace, so it must be $M$. Thus $M$ itself is an infinite-dimensional minimal non-trivial $U$-invariant subspace, contradicting $(iii)$.
\end{proof}

Proposition~\ref{prop:ISP} is a standard motivation for studying concrete universal operators:
a positive solution of the Hilbert space invariant subspace problem would follow if one could
show that the minimal invariant subspaces of \emph{one} universal operator (hence of any universal
operator) are $1$-dimensional; see \cite{CP2011,CGGintro}.

\begin{remark}\label{obstruction}
    In the Banach space setting, things are rather different. Some implications of Proposition \ref{prop:ISP} remain true: $(ii) \iff (iii)$ and $(iii) \implies (i)$ since its proof uses no Hilbert-space specific feature; this is not the case for $(i) \implies (iii)$. Indeed, we would need that every infinite-dimensional minimal non-trivial $U$-invariant subspace is isomorphic to $X$, and we do not know any such examples. In addition, it is not clear that such a hypothesis on $X$ does not force $X$ to be isometric to a Hilbert space... 

    Let us discuss this obstruction a little further. Consider the following property for $U$ universal on a Banach space $X$:
    \begin{center}
    $(H_U)$ ``Every infinite-dimensional minimal non-trivial $U$-invariant subspace is isomorphic to $X$''
\end{center}
Now, the question is: 

\begin{question}\label{questionHU}
  If $(H_U)$ is true for one Banach space $X$ and one universal operator $U$, do we have $X$ isomorphic to a Hilbert space?  
\end{question} Note that the answer is yes if one removes the word ``minimal". Indeed, let $F \subset X$ be a closed infinite-dimensional subspace. Using the definition of $U$ universal with $T=0$, we obtain $J : X \longrightarrow \ker(U)$ an isomorphic embedding. Then, $J(F)$ is a non-trivial $U$-invariant subspace ($J(F) \subset \ker(U) \neq X$ since $U \neq 0$), so $F \simeq J(F) \simeq X$. Hence, $X$ is isomorphic to all its infinite-dimensional closed subspaces, and by the solution of the homogeneous space problem \cite{Gowers}, \cite{KT}, $X$ is isomorphic to $\ell_2$.

Note also that a positive answer to Question \ref{questionHU} would give counterexamples to the invariant subspace problem in the reflexive non-isomorphic to Hilbert space setting. Indeed, let $X$ be reflexive but non-isomorphic to a Hilbert space and admitting a universal operator $U$ (for instance $\ell_p, p \neq 2$). If every operator on every infinite-dimensional closed subspace of $X$ has a non-trivial invariant subspace, there would be no infinite-dimensional minimal $U$-invariant subspace, so $(H_U)$ would hold vacuously, and a positive answer to Question \ref{questionHU} would give $X \simeq \ell_2$, a contradiction.
    
\end{remark}
    
In Section \ref{univBanach}, we will define another notion of universality, more compatible with the Banach space setting, and obtain a result similar to Proposition \ref{prop:ISP} in that setting.

Nevertheless, for a universal operator in the Banach space setting, one has the following proposition:
 
\begin{proposition}
    Let $X$ be an infinite-dimensional complex Banach space and $U$ universal on $X$. Then the following properties are equivalent:
    \begin{enumerate}
        \item Every operator on $X$ has a non-trivial invariant subspace.
        \item There is no minimal non-trivial $U$-invariant subspace isomorphic to $X$.
    \end{enumerate}
\end{proposition}

\begin{proof}
    $(1) \implies (2)$: Assume that there exists a minimal non-trivial $U$-invariant subspace isomorphic to $X$. Denote it $M$. Let $J : X \longrightarrow M$ be an isomorphism. Then, $T:=J^{-1} \circ U\restriction_{M} \circ J \in \mathcal{L}(X)$ does not admit any non-trivial invariant subspace, otherwise, it would contradict the minimality of $M$.
    \smallskip

    $(2) \implies (1)$: Assume that there exists an operator $T \in \mathcal{L}(X)$ with no non-trivial invariant subspace. Then, since $U$ is universal, there exists a closed subspace $M \subset X$ which is $U$-invariant, a complex constant $c \neq 0$, and an isomorphism $J : M \longrightarrow X$ such that:
    $$U\restriction_{M}=J^{-1}(cT)J.$$
    Again, the subspace $M$ is nonzero and proper. Every nonzero closed invariant subspace of $M$ pulls back to a nonzero closed $T$-invariant subspace, so the only closed $U$-invariant subspaces contained in $M$ are $\{0\}$ and $M$; that is, $M$ is a minimal non-trivial $U$-invariant subspace. This contradicts $(ii)$ since $M$ is isomorphic to $X$.
\end{proof}

\subsection{Caradus' criterion on Hilbert space and examples}\label{Caradusetex}

In this subsection, for pedagogical purposes, we recall the famous Caradus criterion of universality in the Hilbert frame, and some interesting examples of universal operators that have already been provided in the past.

\begin{theorem}[Caradus \cite{Caradus1969}]\label{thm:caradus-hilbert}
Let $H$ be a separable infinite-dimensional Hilbert space and let $U\in\mathcal L(H)$.
Assume that $U$ is surjective and $\dim\ker(U)=\infty$. Then $U$ is universal on $H$.
\end{theorem}

\begin{remark}
On a Hilbert space, every closed subspace is complemented by its orthogonal complement.
In particular, $\ker(U)$ is complemented. If $U$ is surjective, then $U$ admits a bounded
right inverse: indeed, the restriction $U\restriction_{\ker(U)^\perp}:\ker(U)^\perp\to H$
is a bounded bijection, hence has a bounded inverse by the open mapping theorem.
Thus Theorem~\ref{thm:caradus-hilbert} is a Hilbert space instance of the Banach space theorem
proved later (Theorem~\ref{thm:caradus-banach}).
\end{remark}

We turn to a first classical example: a shift of infinite multiplicity. This example will be used repeatedly in the sequel.

\begin{example}[Rota's model operator \cite{Rota1959, Rota1960}]\label{ex:shift}
Let
$$
H = \Big(\bigoplus_{n\geq 0}\ell_2\Big)_{\ell_2}
= \Big\{x=(x_n)_{n\geq 0}: x_n\in\ell_2,\ \sum_{n\geq 0}\norm{x_n}^2<\infty\Big\},
$$
equipped with $\norm{x}^2=\sum_{n\geq 0}\norm{x_n}^2$. Define the \emph{backward shift of
infinite multiplicity} $B_\infty\in\mathcal L(H)$ by
$$
B_\infty(x_0,x_1,x_2,\dots) = (x_1,x_2,\dots).
$$
Then $B_\infty$ is surjective and $\ker(B_\infty)=\{(x_0,0,0,\dots):x_0\in\ell_2\}$ is
infinite-dimensional. Hence $B_\infty$ is universal on $H$ by Theorem~\ref{thm:caradus-hilbert}.
\end{example}

A similar example is the following translation operator on 
$L^2(0,\infty)$. This operator is actually unitary equivalent to the backward shift above, via the unitary map $f \longmapsto (f(n+\cdot)\restriction_{[0,1]})_{n \geq 0}$ from $L^2(0,\infty)$ onto $\ell_2(L^2(0,1)) \simeq H$. We still mention it for convenience:

\begin{example}
Following \cite[Section~3.2.1]{CP2011}, one may consider the translation operator $S_1$ on $L^2(0,\infty)$ defined by
$$
(S_1f)(t)=f(t+1),\qquad t\geq 0.
$$
It is surjective (given $g$, take $f(t)=0$ for $t\in[0,1)$ and $f(t)=g(t-1)$ for $t\geq 1$),
and $\ker(S_1)$ contains all functions supported in $[0,1]$, hence is infinite-dimensional.
Thus $S_1$ is universal as well.    
\end{example}

Next, let $\mathbb D$ be the unit disc and let $\mathbb H^2(\mathbb D)$ denote the Hardy-Hilbert space.
For an analytic self-map $\varphi:\mathbb D\to\mathbb D$, the \emph{composition operator}
$C_\varphi$ on $\mathbb H^2(\mathbb D)$ is given by $C_\varphi f = f\circ \varphi$.

\begin{example}[Nordgren--Rosenthal--Wintrobe \cite{NRW1987}]\label{ex:composition}
Let $\varphi$ be the hyperbolic automorphism
$$
\varphi(z)=\frac{z+\frac12}{1+\frac12 z},\qquad z\in\mathbb D,
$$
and consider the operator
$$
U = C_\varphi - I\in\mathcal L\big(\mathbb H^2(\mathbb D)\big).
$$
Then $U$ is universal on $\mathbb H^2(\mathbb D)$.

Indeed, Nordgren, Rosenthal and Wintrobe proved that if $\varphi$ is a hyperbolic automorphism,
then $C_\varphi-\lambda I$ is universal on $\mathbb H^2(\mathbb D)$ whenever
$\lambda$ lies in the interior of the spectrum of $C_\varphi$ \cite{NRW1987} (see also \cite{CGGnewproof} for a new proof).
Here, $\varphi(\pm 1)= \pm 1$ and $\varphi'(1)=\frac{1}{3}$, so $\sigma(C_{\varphi})=\{ \lambda \in \mathbb{C} : \frac{1}{\sqrt{3}} \leq |\lambda | \leq \sqrt{3}\}$ and $1$ lies in its interior, so $U=C_\varphi-I$
is universal.
\end{example}

\section{Hereditarily universal operators}\label{univBanach}
In this part, we define a new stronger notion of universality, called \emph{hereditary universality}, in order to generalize Proposition \ref{prop:ISP}. We also show that the Banach space version of Caradus' criterion, which was already known, actually provides good sufficient conditions for an operator to be hereditarily universal. Then, we discuss some geometrical properties of infinite-dimensional Banach spaces which guarantee the hereditary universality of an operator which is surjective and whose kernel is infinite-dimensional and complemented.

\subsection{Definition and criterion}

\begin{definition}
    Let $X$ be a complex Banach space. We say that $U \in \mathcal{L}(X)$ is \emph{hereditarily universal} if for every infinite-dimensional closed subspace $Y \subset X$ and every $T \in \mathcal L(Y)$, there exist $c\in \C\setminus \{0\}$, a closed $U$-invariant subspace $M \subset X$, and an isomorphism $W : Y \to M$ such that $U\restriction_{M}W = cWT$. 
\end{definition}

It is clear that a hereditarily universal operator on a Banach space is in particular universal in the sense of Rota (Definition \ref{defuniv}). The following theorem gives sufficient conditions for an operator to be hereditarily universal. It is largely inspired by the Banach space version of Caradus' proof (compare \cite{Caradus1969,CP2011}).

\begin{theorem}\label{thm:caradus-banach}
Let $X$ be an infinite-dimensional Banach space and let $U\in\mathcal L(X)$.
Assume that:
\begin{enumerate}[(a)]
\item $U$ admits a bounded right inverse, i.e.\ there exists $V\in\mathcal L(X)$ such that $UV=I_X$;
\item $\ker(U)$ contains a closed subspace isomorphic to $X$.
\end{enumerate}
Then $U$ is hereditarily universal on $X$.
\end{theorem}

\begin{proof}
    Fix $V\in\mathcal L(X)$ with $UV=I_X$.
    By (b), there exists a closed subspace $K\subset\ker(U)$
    and a bounded linear isomorphism $J:X\to K$ with
    bounded inverse $J^{-1}:K\to X$.

    Let $Y \subset X$ be a closed subspace and let $T\in\mathcal L(Y)$ be arbitrary. Choose a scalar
    $c\in\C\setminus\{0\}$ so that
    $$
    |c|\|V\|\|T\|\|J\|\|J^{-1}\|<\frac12.
    $$
    Put $t=|c|\|V\|\|T\|$. Since
    $\|J\|\|J^{-1}\|\geq1$, we have $t<\frac12$. Define $A:Y \longrightarrow X$ by the norm-convergent series
    $$
    Ay=\sum_{n=0}^{\infty}c^nV^nJT^ny,\qquad y \in Y.
    $$
    Indeed, $\|c^nV^nJT^n\|\leq t^n\|J\|$, so the series converges in operator norm and defines
    a bounded operator with
    $$
    \|A\|\leq\frac{\|J\|}{1-t}.
    $$

    \smallskip
    Recall that $UV=I_X$, while $UJ=0$ since
    $J(X)\subset\ker(U)$. Moreover, the series gives
    $A=J+cVAT$. Hence
    $$
    UA=UJ+cUVAT=A(cT).
    $$

    Let $y \in Y$. The triangle inequality and the fact
    that $J$ has bounded inverse yield
    $$
    \begin{aligned}
        \|Ay\|
        &\geq\|Jy\|-|c|\|VATy\|\\
        &\geq
        \bigl(\beta-|c|\|V\|\|A\|\|T\|\bigr)\|y\|,
    \end{aligned}
    $$
    where $\beta=\frac{1}{\|J^{-1}\|}$. Moreover,
    $$
    |c|\|V\|\|A\|\|T\|
    \leq\frac{t\|J\|}{1-t}<\beta,
    $$
    since $t\|J\|\|J^{-1}\|<\frac12$ and
    $1-t>\frac12$. Consequently,
    $$
    \|Ay\|\geq\alpha\|y\|,
    \qquad
    \alpha=\beta-\frac{t\|J\|}{1-t}>0.
    $$
    Thus $A$ is bounded below. In particular,
    $\mathrm{Ran}(A)$ is closed, and $A$ is a bounded
    isomorphism onto its range.

    \smallskip
    \noindent\emph{Conclusion.}
    Let $M=\mathrm{Ran}(A)$, a closed subspace of $X$.
    Since $UA=A(cT)$ and $A(Y)=M$, we have
    $U(M)\subset M$, so $M$ is $U$-invariant.
    The intertwining identity can be rewritten as
    $$
    U\restriction_M\,A=A(cT),
    $$
    and since $A:Y \longrightarrow M$ is a bounded isomorphism,
    it follows that $U\restriction_M$ is similar to $cT$.
    Since $T$ and $Y$ were arbitrary, $U$ is hereditarily universal.
\end{proof}

We can now state an analogue of Proposition \ref{prop:ISP} in the Banach space setting:

\begin{proposition}\label{ISPBanach} Let $X$ be a complex Banach space and $U \in \mathcal{L}(X)$ be hereditarily universal. Then, the following properties are equivalent:
\begin{enumerate}
    \item[(i)] For every infinite-dimensional closed subspace $Y \subset X$, every $T \in \mathcal{L}(Y)$ has a non-trivial invariant subspace.
    \item[(ii)] every minimal non-trivial $U$-invariant subspace has finite dimension;
    \item[(iii)] every minimal non-trivial $U$-invariant subspace has dimension $1$.
\end{enumerate}
    
\end{proposition}

\begin{proof}
We already showed that $(ii) \iff (iii)$. For $(iii) \implies (i)$, we proceed the same way as to show $(iii) \implies (i)$ in Proposition \ref{prop:ISP} but we use the hereditary universality of $U$ instead of universality.

Let us show that $(i) \implies (ii)$. Let $M$ be a minimal non-trivial $U$-invariant subspace of $X$. Assume that $M$ is infinite-dimensional. Then, by $(i)$,  $U\restriction_{M} \in \mathcal{L}(M)$ has a non-trivial invariant subspace, but it contradicts the minimality of $M$. 
\end{proof}

\begin{remark}
   As in the Hilbert case, if one could find a hereditarily universal operator on a Banach space $X$ whose minimal non-trivial invariant subspaces all have dimension 1, Proposition \ref{ISPBanach} would give a positive answer to the invariant subspace problem for all operators on all infinite-dimensional closed subspace of $X$, for instance for $X=\ell_p$, $1<p<\infty$, $p \neq 2$, a case which is still open. A single separable space cannot settle the problem for all reflexive spaces at once. Indeed, no separable reflexive space contains isomorphic copies of all separable reflexive space \cite{Szlenk}. Moreover, a separable space containing such copies contains a copy of $\ell_1$ \cite{Bourgain}, hence, by Read's counterexample \cite{Read2}, an infinite-dimensional closed subspace carrying an operator without non-trivial invariant subspace; on such a space, condition $(i)$ of Proposition \ref{ISPBanach} fails, so such an operator cannot exist. 
   \smallskip 

    Note also that in the Hilbert case, a universal operator on $H$ is in fact hereditarily universal since every infinite-dimensional closed subspace of a separable infinite-dimensional Hilbert space $H$ is isometrically isomorphic to $H$.

    \end{remark}

\subsection{Large complemented subspaces and the Caradus criterion}

We now discuss some classes of Banach spaces for which the assumptions of
Theorem~\ref{thm:caradus-banach} can be verified in a natural way. We begin by
recording a simple observation which reformulates the existence of a bounded
right inverse.

\begin{lemma}\label{lem:rightinverse}
Let $X$ be a Banach space and let $U\in\mathcal L(X)$ be surjective. The following
are equivalent:
\begin{enumerate}[(i)]
\item $U$ admits a bounded right inverse;
\item $\ker(U)$ is complemented in $X$.
\end{enumerate}
\end{lemma}

\begin{proof}
Suppose first that $V\in\mathcal L(X)$ satisfies $UV=I_X$. Then
$P=I_X-VU$ is a bounded projection onto $\ker(U)$.

Conversely, assume that $\ker(U)$ is complemented and write
$X=\ker(U)\oplus Z$ for some closed subspace $Z$ of $X$. Since $U$ is
surjective, the restriction $U\restriction_Z:Z\longrightarrow X$
is a bounded bijection. By the open mapping theorem, its inverse is bounded.
Thus
$$
V=(U\restriction_Z)^{-1}:X\longrightarrow Z\subset X
$$
is a bounded right inverse of $U$.
\end{proof}

In view of this lemma, for a surjective operator the hypothesis of Theorem~\ref{thm:caradus-banach} reduce to the requirement that its kernel is complemented and contains a copy of the
ambient space. This leads naturally to the following property:
\begin{equation}\label{eq:large-complemented}
\begin{split}
&\text{every infinite-dimensional complemented subspace of $X$}\\
&\text{contains a closed subspace isomorphic to $X$.}
\end{split}
\tag{$\mathcal C$}
\end{equation}

\begin{proposition}\label{prop:large-kernel-caradus}
Let $X$ be an infinite-dimensional Banach space satisfying
\eqref{eq:large-complemented}, and let $U\in\mathcal L(X)$ be surjective.
If $\ker(U)$ is infinite-dimensional and complemented in $X$, then $U$ is hereditarily
universal on $X$.
\end{proposition}

\begin{proof}
Since $U$ is surjective and $\ker(U)$ is complemented,
Lemma~\ref{lem:rightinverse} shows that $U$ admits a bounded right inverse.
Moreover, $\ker(U)$ is an infinite-dimensional complemented subspace of $X$,
so \eqref{eq:large-complemented} yields a closed subspace $Z\subset\ker(U)$ which is isomorphic to $X$, and $U$ is hereditarily universal by
Theorem~\ref{thm:caradus-banach}.
\end{proof}

Several familiar structural properties of Banach spaces imply
\eqref{eq:large-complemented}. Recall that an infinite-dimensional Banach
space $X$ is called \emph{minimal} if every infinite-dimensional closed
subspace of $X$ contains a further closed subspace isomorphic to $X$. It is
called \emph{complementably minimal} if this copy of $X$ can moreover be
chosen complemented in $X$. Finally, $X$ is called \emph{prime} if every
infinite-dimensional complemented subspace of $X$ is isomorphic to $X$.

Both minimality and primeness imply \eqref{eq:large-complemented}, while
complementable minimality implies minimality. Notice, however, that minimality
and primeness express somewhat different properties. Minimality concerns all
infinite-dimensional subspaces, but only requires them to contain a copy of
$X$, whereas primeness concerns only complemented subspaces, but requires
them to be isomorphic to $X$.

We therefore obtain the following immediate consequence.

\begin{corollary}\label{cor:prime-minimal-caradus}
Let $X$ be an infinite-dimensional Banach space which is either minimal or
prime, and let $U\in\mathcal L(X)$ be surjective. If $\ker(U)$ is
infinite-dimensional and complemented in $X$, then $U$ is hereditarily universal on $X$.
\end{corollary}

In particular, the conclusion holds whenever $X$ is complementably minimal.

\subsubsection{Classical sequence spaces}

The classical sequence spaces provide the most immediate examples. It is
well known that $c_0$ and $\ell_p$, $1\leq p<\infty$, are complementably
minimal; see, for instance, \cite[Chapter~2]{LT1977}. In particular, they are
minimal. Moreover, a classical theorem of Lindenstrauss asserts that
$\ell_\infty$ is prime; see \cite{Lindenstrauss1967}. Hence we obtain the
following.

\begin{corollary}\label{cor:classical-sequence}
Let
$$
X=c_0,\qquad X=\ell_p\quad (1\leq p<\infty),
\qquad\text{or}\qquad X=\ell_\infty,
$$
and let $U\in\mathcal L(X)$ be surjective. If $\ker(U)$ is
infinite-dimensional and complemented in $X$, then $U$ is hereditarily universal on $X$.
\end{corollary}

\begin{proof}
For $c_0$ and $\ell_p$, this follows directly from
Proposition~\ref{prop:large-kernel-caradus}. If $X=\ell_\infty$, the kernel is an
infinite-dimensional complemented subspace of $\ell_\infty$, hence it is
isomorphic to $\ell_\infty$ by primeness. Lemma~\ref{lem:rightinverse} and
Theorem~\ref{thm:caradus-banach} then yield the conclusion.
\end{proof}

\begin{example}[A simple hereditarily universal operator on $\ell_p$]\label{ex:lp-simple}
Fix $1\leq p<\infty$. Define $U\in\mathcal L(\ell_p)$ on the canonical basis
$(e_n)_{n\geq 0}$ by
$$
U(e_{2n})=0,\qquad U(e_{2n+1})=e_n,\qquad n\geq 0.
$$
Equivalently,
$$
U(x_0,x_1,x_2,x_3,\dots)=(x_1,x_3,x_5,\dots).
$$
The operator $U$ is surjective, with bounded right inverse
$$
V(y_0,y_1,\dots)=(0,y_0,0,y_1,0,y_2,\dots),
$$
and
$$
\ker(U)=\{(x_0,0,x_2,0,x_4,0,\dots)\}
$$
is an infinite-dimensional complemented subspace of $\ell_p$. Therefore $U$
is hereditarily universal by Corollary~\ref{cor:classical-sequence}.

This example also has the property that all its powers are hereditarily universal. Indeed,
for every $k\geq 1$,
$$
U^k(x)
=
\big(x_{2^k-1},x_{2^{k+1}-1},x_{3\cdot 2^k-1},\dots\big),
$$
and hence
$$
\ker(U^k)
=
\big\{x\in\ell_p:
x_{2^kn+2^k-1}=0\text{ for every }n\geq0\big\}.
$$
Thus $\ker(U^k)$ is again an infinite-dimensional complemented coordinate
subspace of $\ell_p$, and Corollary~\ref{cor:classical-sequence} shows that
$U^k$ is hereditarily universal for every $k\geq1$.

The same coordinate construction works, without any change, on $c_0$.
\end{example}

Let us finally mention that the above scheme applies well beyond the classical sequence spaces. For instance, Schlumprecht's space $S$ is complementably minimal \cite{AndroulakisSchlumprecht2003}, while the dual $T^*$ of Tsirelson's space is minimal \cite{CasazzaJohnsonTzafriri1984}; hence Proposition~\ref{prop:large-kernel-caradus} applies to both spaces.

For $1\leq p<\infty$, $p\neq2$, the situation is different on $L^p[0,1]$: complemented copies of $\ell_2$ (when $p>1$) or $\ell_1$ (when $p=1$) contain no copy of $L^p[0,1]$. Thus, infinite-dimensionality and complementedness alone do not imply the kernel condition required by $(\mathcal C)$, namely $L^p[0,1]\hookrightarrow\ker(U)$. This embedding is also necessary for universality, as follows by modelling the zero operator. However, these subspace examples alone do not construct a surjective operator with such a kernel and therefore do not provide counterexamples to universality. Nevertheless, since $L^p([0,1]) \simeq \ell_p(L^p([0,1]))$, the backward shift on $\ell_p(L^p([0,1]))$ satisfies the hypothesis of Theorem \ref{thm:caradus-banach}, hence $L^p([0,1])$ carries hereditarily universal operators. More generally, this applies to every $X$ with $X \simeq \ell_p(X)$ or $X \simeq c_0(X)$.

\subsection{A structural obstruction: spaces with no universal operators}

Universality forces $\ker(U)$ to contain a copy of $X$: indeed, take $T=0$ in the definition.
Thus, if a space has the property that every non-compact operator has ``small kernel'' (for instance,
finite-dimensional kernel), then it cannot support a universal operator.

A striking example is the Argyros--Haydon scalar-plus-compact space $X_{AH}$ \cite{AH2011}, where
every $T\in\mathcal L(X_{AH})$ has the form $T=\lambda I + K$ with $K$ compact. If $\lambda\neq 0$,
then $\lambda I+K$ is Fredholm (hence has finite-dimensional kernel), while if $\lambda=0$ then $T$
is compact and cannot be universal (a universal operator $U$ models $I_{X_{AH}}$, so $U\restriction_{M}$ is similar to $cI_{X_{AH}}$ for some infinite-dimensional closed $M$; thus $U$ is bounded below on a infinite-dimensional subspace and is neither compact nor strictly singular). Consequently $X_{AH}$ admits no universal operator.

\begin{remark}
In the same space, let $\mathcal{A}=\{\lambda I+K : \lambda \in \mathbb{C}^*, K \text{ compact operator}\}$. Then, there is no universal operator for the class $\mathcal{A}$, that is, no operator $U \in \mathcal{L}(X_{AH})$ such that for every $T \in \mathcal{A}$, there are $c \neq 0$ and a closed $U$-invariant subspace $M$ with $U\restriction_{M}$ similar to $cT$. Indeed, assume there exists $U$ universal for $\mathcal{A}$. Let $U=\lambda_0I+K_0$. Since $I \in \mathcal{A}$, $U\restriction_{M}$ is similar to $cI$ for some infinite-dimensional closed subspace $M \subset X_{AH}$, so $U$ is not compact, hence $\lambda_0 \neq 0$ and $U$ is Fredholm. We denote $n_0=\dim(\ker(\lambda_0I+K_0))$. Then, for every closed $U$-invariant subspace $M \subset X_{AH}$, $\dim(\ker(U\restriction_\mathcal{M})) \leq n_0$.
Let $(x_1,...,x_{n_0+1})$ be a family of linearly independent vectors, and let $G$ be a closed subspace of $X_{AH}$ such that $X_{AH}=\mathrm{span}(x_1,...,x_{n_0+1}) \oplus G$. Let us define an operator $K_1$ such that: 
$$K_1\restriction_G=0 \text{ and } \forall j \in \{1,...,n_0+1\}, K_1(x_j)=-\lambda_0x_j$$
Then $K_1$ is compact since it has finite rank. Let $T=\lambda_0 I+K_1$. Then, $T \in \mathcal {A}$ and $\dim(\ker(T))=n_0+1>n_0$, so no restriction of $U$ can be similar to $\mu T$ for $\mu \in \mathbb{C}^*$, a contradiction.
\end{remark}

\begin{remark}
The article \cite{GM} of Gowers and Maurey presents a construction of a separable Banach space $X$ which contains no unconditional basic sequence. The authors also show that every bounded linear operator $T$ from $X$ into $X$ can be written $T=\lambda I+S$, where $\lambda \in \mathbb{C}$ and $S$ is strictly singular. As a consequence, there cannot exist any universal operator in $X$. Indeed, if $\lambda \neq 0$, we have that $\ker(\lambda I)=\{0\}$ has finite dimension and $\mathrm{Ran}(\lambda I)=X$ has finite codimension, so $\lambda I$ is Fredholm. We know that the sum of a Fredholm operator and a strictly singular operator is Fredholm (see \cite{LT1977}), so $T$ is Fredholm and cannot be universal. If $\lambda=0$, $T=S$ is strictly singular and cannot be universal.
\end{remark}

\section{Spectrum of a universal operator}\label{spectrum}
In this section, we study spectral properties of universal operators (in the sense of Rota, Definition \ref{defuniv}) on a Banach space. It is inspired by the work of Rikka Schröderus and Hans-Olav Tylli in \cite{schro}.
For a Banach space $X$, we will use the same notation as in \cite{schro}: Let $\mathcal{U}(X)=\{U \in \mathcal{L}(X) : U \text{ is universal}\}$, and:
$$\Phi_-(X)=\{S \in \mathcal{L}(X) : \dim(X/\mathrm{Ran}(S))<\infty, \;  \mathrm{Ran}(S) \text{ is closed}\};$$
$$\Phi_+(X)=\{S \in \mathcal{L}(X) : \dim(\ker(S))<\infty,\; \mathrm{Ran}(S) \text{ is closed}\}.$$
The set $\Phi(X):=\Phi_-(X) \cap \Phi_+(X)$ represents the set of all the Fredholm operators. We also define the $\sigma_e^+$-spectrum of $S \in \mathcal{L}(X)$ by:
$$\sigma_e^+(S,X)=\{\lambda \in \mathbb{C} ~  | ~  S-\lambda I \notin \Phi_+(X)\}.$$
It is known that $\sigma_e^+(S,X)$ is a nonempty compact subset of the essential spectrum:
$$\sigma_e(S,X)=\{\lambda \in \mathbb{C} ~  | ~ S-\lambda I \notin \Phi(X)\}.$$
Finally, let us denote $\sigma_p(S,X)$ the point spectrum of $S \in \mathcal{L}(X)$.

\begin{theorem}
Let $1 \leq p<\infty$. Let $U \in \mathcal{U}(\ell_p)$. Then:
\begin{enumerate}
    \item[(i)] There exists $r>0$ such that 
    $$B(0,r) \subset \sigma_p(U,\ell_p) \cap \sigma_e^+(U,\ell_p) \subset \sigma_e(U,\ell_p) \subset \sigma(U,\ell_p)$$
    Moreover, any $\lambda \in B(0,r)$ is an eigenvalue of $U$ having infinite multiplicity. In particular, $0$ is an interior point of the subsets $\sigma_e^+(U,\ell_p), \sigma_e(U,\ell_p),\\\sigma_p(U,\ell_p)$ and $\sigma(U,\ell_p)$.
    \item[(ii)] There exist $r>0$ and a holomorphic vector-valued map $z \in B(0,r) \longmapsto y_z \in \ell_p\setminus \{0\}$ such that:
    $$\forall z \in B(0,r), U(y_z)=zy_z$$
\end{enumerate}
\end{theorem}

\begin{proof}
Let $q=\frac{p}{p-1}$ be the conjugate exponent of $p$, with the convention $q=\infty$ if $p=1$.
Let us denote 
$$
X_0=\left(\bigoplus\limits_{\mathbb{N} \cup \{0\}} \ell_p\right)_{\ell_p}=\ell_p(\ell_p).
$$ 
We know that $\ell_p \simeq X_0$ (see \cite{AK}). Let us consider the backward shift on $X_0$:
$$\fonction{B_{\infty}}{X_0}{X_0}{x=(x_0,x_1, x_2,...)}{(x_1,x_2,...).}$$
Then, for $\lambda \in \mathbb{D}$, $x \in X_0, x \neq 0$:
$$B_{\infty}(x)-\lambda x=(x_1-\lambda x_0,x_2-\lambda x_1,...),$$
hence, $\ker(B_{\infty}-\lambda I)=\{(x_0,\lambda x_0,\lambda^2x_0,...) : x_0 \in \ell_p\}$. This shows that $\lambda$ is an eigenvalue of $B_{\infty}$ of infinite multiplicity, so $\mathbb{D} \subset \sigma_p(B_{\infty},X_0) \cap \sigma_e^+(B_{\infty},X_0)$. We also  have that:
$$\ker(B_{\infty})=\{(x,0,0,...) : x \in \ell_p\} \simeq \ell_p.$$
Since $X_0 \simeq \ell_p$ and $U$ is universal, there exists a closed $U$-invariant subspace $M \subset \ell_p$, a constant $c \neq 0$ and an isomorphism $J : X_0 \longrightarrow M$ such that $U\restriction_M=J(cB_{\infty})J^{-1}$ (universality being applied to the operator on $\ell_p$ that corresponds to $B_{\infty}$ via $X_0 \simeq \ell_p$. Since all the spectral elements are preserved by isomorphism:
$$B(0,|c|) \subset  \sigma_p(cB_{\infty},X_0)=\sigma_p(U\restriction_M,M) \subset \sigma_{p}(U,\ell_p),$$
and
$$B(0,|c|) \subset \sigma_e^+(cB_{\infty},X_0)=\sigma_e^+(U\restriction_M,M) \subset \sigma_e^+(U,\ell_p),$$
where the last inclusion holds because, for $\lambda \in B(0,|c|)$, $J$ maps the infinite-dimensional space $\ker(cB_{\infty}-\lambda I)$ into $\ker(U-\lambda I)$, hence $U-\lambda I \notin \Phi_+(\ell_p)$ and $\lambda$ is an eigenvalue of $U$ of infinite multiplicity. This proves $(i)$ with $r=|c|$.

We turn to the proof of $(ii)$.  For $z \in \mathbb{D}$ and $x_0 \in \ell_p, x_0 \neq 0$, let $x_z:=(z^nx_0)_n \in X_0$. The map $z \longmapsto x_z$ is weakly holomorphic, since:
$$\forall y=(y_n) \in \left(\bigoplus_{\mathbb{N} \cup \{0\}} \ell_q\right)_{\ell_q}, z \longmapsto \sum_{n=0}^{+\infty} z^n\langle x_0,y_n\rangle_{\ell_p,\ell_q} \text{ is analytic.}$$
Hence, it is holomorphic (see \cite[Theorem 8.12]{muj}). 

Now for $z\in B(0,|c|)$, set $y_z=Jx_{z/c}$. Then the map $z\mapsto y_z$ is holomorphic, and
$$
U(y_z) = cJB_\infty x_{z/c} = c\frac{z}{c}Jx_{z/c} = zy_z.
$$
This proves (ii).
\end{proof}

\begin{remark}
Since $c_0(c_0)$ is also isomorphic to $c_0$, the same argument applies to a universal operator $U$ on $c_0$.
\end{remark}

\section{Universal \texorpdfstring{$n$}{n}-tuples of operators}

This is based on Vladimir Müller's article \cite{Mul} on the same subject. His work deals with operators on Hilbert spaces. We present here a generalization to Banach spaces. In the following, $X$ and $Y$ denote infinite-dimensional complex Banach spaces.

\subsection{Commuting \texorpdfstring{$n$}{n}-tuples}

\begin{definition}\label{univtuples}
Let $X,Y$ be fixed infinite-dimensional complex Banach spaces. A commuting tuple $T=(T_1,...,T_n)\in \mathcal L(Y)^n$ is universal for commuting tuples on $X$ if, for every commuting $A=(A_1,...,A_n)\in \mathcal L(X)^n$, there exist $c\neq0$, a closed common invariant subspace $M\subset Y$, and an isomorphism $J:X\to M$ such that $T_j\restriction_{M}J=cJA_j$ for all $j$.
\end{definition}

Here, the $n$-tuple $T=(T_1,\ldots,T_n)$ is viewed as a family of commuting operators acting simultaneously on the same Banach space $Y$, rather than as a single operator acting on the product space $Y^n$.
In the light of Theorem 3 of \cite{Mul}, the following theorem gives sufficient conditions for an $n$-tuple $(T_1,...,T_n)$ to be universal.

\begin{theorem}\label{tuples}
Let $X,Y$ be complex infinite-dimensional Banach spaces. Let $n \geq 1$ and let $T_1,...,T_n \in \mathcal L(Y)$ be a commuting $n$-tuple of operators satisfying 
\begin{enumerate}
    \item[(i)] $\dim\left(\displaystyle \bigcap_{j=1}^{n} \ker(T_j)\right)=\infty;$
    \item[(ii)] for all $F \subset \{1,...,n\}$ and $j \in \{1,...,n\} \setminus F$, the restriction $T_j\restriction_{\bigcap_{i \in F}\ker(T_i)}$ is surjective, i.e., $T_j\left(\displaystyle\bigcap_{i \in F} \ker(T_i)\right)=\displaystyle\bigcap_{i \in F} \ker(T_i)$. Note that for $F=\emptyset$, it means that $T_j$ is surjective;
    \item[(iii)] $\forall F \subset \{1,...,n\}$ and $j \in \{1,...,n\} \setminus F$, $\displaystyle\bigcap_{i \in F} \ker(T_i) \cap \ker(T_j)$ is complemented in $Y$. Again, if $F=\emptyset$, it means that $\ker(T_j)$ is complemented in $Y$;
    \item[(iv)] There exists an operator denoted $V_{(0,...,0)} : X \longrightarrow \displaystyle\bigcap_{i=1}^{n} \ker(T_i)$ such that $\|V_{(0,...,0)}(x)\| \geq \|x\|$ for all $x \in X$. 
\end{enumerate}
Then, $(T_1,...,T_n)$ is universal for commuting $n$-tuples of operators on $X$.
\end{theorem}

\begin{remark}
      Note that $(iv) \implies (i)$. Indeed, the operator $V_{(0,...,0)}$ is in particular injective, so $+\infty=\dim(X) \leq \dim\left(\displaystyle\bigcap_{j=1}^{n} \ker(T_j)\right)$.

      Still, we put the condition $(i)$ to emphasize it because it is also necessary: just apply Definition \ref{univtuples} with $(A_1,...,A_n)=(0,...,0)$. Applying Definition \ref{univtuples} to $(0,\ldots,0)$ also shows that the condition $(iv)$ is necessary (take $V_{(0,\ldots,0)}=\|J^{-1}\|J)$.  It would be nice to know in what setting the conditions $(ii)$ and $(iii)$ are necessary.
\end{remark}

The proof requires several lemmas, which are adapted from Lemma 4, Lemma 5 and Lemma 6 of \cite{Mul}. Let $T_1,...,T_n \in \mathcal L(Y)$ be a commuting $n$-tuple satisfying the conditions of Theorem \ref{tuples}. By hypothesis $(ii)$ of Theorem \ref{tuples}, one has the surjectivity of $T_j$.
By $(iii)$ and the open mapping theorem, each such operator admits a bounded right inverse which will be denoted $\widehat{T_j}$ (see Lemma \ref{lem:rightinverse}). So we have $T_j\widehat{T_j}=\mathrm{Id}_Y$.

\begin{lemma}\label{proj}
Let $F \subset \{1,...,n\}$ and $j \in \{1,...,n\} \setminus F$. Then:
$$\bigcap_{i \in F} \ker(T_iT_j)=\ker(T_j)+\bigcap_{i \in F} \ker(T_i)$$
Moreover, there exists a projection $P_{j,F} : \displaystyle\bigcap_{i \in F} \ker(T_iT_j) \longrightarrow \displaystyle\bigcap_{i \in F} \ker(T_iT_j)$ such that $\mathrm{Ran}(P_{j,F})=\displaystyle\bigcap_{i \in F} \ker(T_i)$ and $\ker(P_{j,F}) \subset \ker(T_j)$.
\end{lemma}

\begin{proof}
The ``$\supset$" part is clear and true for every commuting $n$-tuple of operators. 

Let $x \in \displaystyle\bigcap_{i \in F} \ker(T_iT_j)$. Then, $T_j(x) \in \displaystyle\bigcap_{i \in F} \ker(T_i)$. By hypothesis (ii) of Theorem \ref{tuples}, there exists $y \in \displaystyle\bigcap_{i \in F} \ker(T_i)$ such that $T_j(y)=T_j(x)$. Thus, $x-y \in \ker(T_j)$ and $x=x-y+y \in \ker(T_j)+ \displaystyle\bigcap_{i \in F} \ker(T_i)$.
\smallskip

Let us prove the second part. We denote $M=\displaystyle\bigcap_{i \in F} \ker(T_i)$ and $L=\ker(T_j)$. Then, $M+L$ is a closed subspace since the first part shows that it is a finite intersection of closed subspaces. Using hypothesis (iii) of Theorem~\ref{tuples}, let $S$ be a topological complement of $M \cap L$. Let $Z=(M+L) \cap S$. Let us show that $M+L=M \oplus (L \cap Z)$. 

Let $x \in M \cap (L \cap Z)$. Then, $x \in M \cap L$ and $x \in S$, so $x=0$. Hence, $M \cap (L \cap Z)=\{0\}$.

Clearly, $M \subset M+(L \cap Z)$. Let $x \in L$. Then $x$ can be uniquely written as 
$$x=y+z \text{ where } y \in M \cap L \text{ and } z \in S$$
So $z \in L \cap Z$ and $x=y+z \in M \oplus (L \cap Z)$. The reverse inclusion is obvious, hence, $M+L=M \oplus (L \cap Z)$.

Finally, there exists a projection $P_{j,F} : M+L \longrightarrow M+L$ which is continuous by the closed graph theorem, such that 
$$
\mathrm{Ran}(P_{j,F})=M=\displaystyle\bigcap_{i \in F} \ker(T_i) \quad  \text{and} \quad \ker(P_{j,F})=L \cap Z \subset L=\ker(T_j).$$
\end{proof}

For $G \subset \{1,...,n\}$, denote $G^c=\{1,...,n\} \setminus G$ and let us introduce:
$$Y_G=\bigcap_{i \notin G} \ker(T_i) \text{ and } R_j^G=P_{j,G^c} \widehat{T_j}\restriction_{Y_G}, ~ j \in G,$$
For $y \in Y_G$, and $i \notin G$, 
$$T_iT_j\widehat{T_j}y=T_iy=0,$$
so $\widehat{T_j}y$ lies in the domain $\displaystyle\bigcap_{i \notin G} \ker(T_iT_j)$ of $P_{j,G^c}$. Then, $R_j^G$ is well-defined. 
Let us also set
\begin{align*}k=\max\Big\{2;&\max\{\|P_{j,F}\| : F \subset \{1,...,n\}, j \in \{1,...,n\} \setminus F\};\\&\max\{\|\widehat{T_j}\| : j \in \{1,...,n\}\}; \max\{\|R_j^G\|: G \subset \{1,...,n\}, j \in G\}\Big\}
\end{align*}

The following lemma is proved in \cite{Mul} in a Hilbert space setting, and the proof carries over verbatim to the Banach space setting, using Lemma \ref{proj}. We include it here for convenience of the reader.

\begin{lemma}\label{opVF}
Let $X$ be a Banach space. Let $G \subset \{1,...,n\}$, $G \neq \emptyset$. Suppose that for all $F \subsetneq G$, there exist bounded linear operators $V_F : X \longrightarrow Y$ satisfying:
$$
T_jV_F=V_{F \setminus \{j\}}, \quad F \subsetneq G, \; j \in F
$$
Then, there exists an operator $V_G : X \longrightarrow Y$ satisfying, for all $j \in G$:
$$T_jV_G=V_{G \setminus \{j\}}$$
and $\|V_G\| \leq (2k^2)^{\mathrm{card}(G)}\max\{\|V_{G \setminus \{j\}}\| : j \in G\}$.
\end{lemma}

\begin{proof}
We prove the statement by induction on $\mathrm{card}(G)$. 

If $\mathrm{card}(G)=1$, $G=\{m\}$, then set $V_{\{m\}}=\widehat{T_m}V_{\emptyset}$. Then: 
$$T_mV_{\{m\}}=T_m\widehat{T_m}V_{\emptyset}=V_{\emptyset} \text{ and } \|V_{\{m\}}\| \leq \|\widehat{T_m}\|\|V_{\emptyset}\| \leq k\|V_{\emptyset}\|.$$

Let $G \subset \{1,...,n\}$, $\mathrm{card}(G) \geq 2$ and suppose that the statement is true for each $\widetilde{G} \subset \{1,...,n\}$ with $1 \leq \mathrm{card}(\widetilde{G})<\mathrm{card}(G)$. 

Fix an $m \in G$ and let $G'=G \setminus \{m\}$. Consider the operators $W_F=V_{F \cup \{m\}}$ $(F \subsetneq G')$. By the induction assumption, there exists an operator $V' : X \longrightarrow Y$ satisfying $T_jV'=W_{G' \setminus \{j\}}$ for all $j \in G'$. Moreover, 
$$\|V'\| \leq (2k^2)^{\mathrm{card}(G')}\max\{\|V_F\| : \mathrm{card}(F)=\mathrm{card}(G)-1\}.$$
Furthermore, let $V''=\widehat{T_m}V_{G'}$. For all $j \in G'$, we have:
$$T_jT_m(V''-V')=T_jV_{G'}-T_mV_{G \setminus \{j\}}=V_{G \setminus \{j,m\}}-V_{G \setminus \{j,m\}}=0.$$

Hence, by Lemma \ref{proj}: 
$$\mathrm{Ran}(V''-V') \subset \bigcap_{j \in G'} \ker(T_jT_m)=\ker(T_m)+\bigcap_{j \in G'} \ker(T_j).$$

Let $P_{m,G'} : \displaystyle\bigcap_{j \in G'} \ker(T_jT_m) \longrightarrow \displaystyle\bigcap_{j \in G'} \ker(T_jT_m)$ the projection given by\\ Lemma~\ref{proj}, satisfying $\mathrm{Ran}(P_{m,G'})=\displaystyle\bigcap_{j \in G'} \ker(T_j)$ and $\ker(P_{m,G'}) \subset \ker(T_m)$. Set
$$V_G=V''+(I-P_{m,G'})(V'-V'')=V'+P_{m,G'}(V''-V').$$
Since $\mathrm{Ran}(I-P_{m,G'})=\ker(P_{m,G'}) \subset \ker(T_m)$, we have $$T_mV_G=T_mV''=V_{G'}=V_{G \setminus \{m\}}.$$
For $j \in G'$, we have $T_jV_G=T_jV'=V_{G \setminus \{j\}}$. Moreover, 
\begin{align*}
\|V_G\| &\leq \|V'\|+\|P_{m,G'}\|(\|V''\|+\|V'\|)\\
&\leq \left((2k^2)^{\mathrm{card}(G)-1}
+k\left(k+(2k^2)^{\mathrm{card}(G)-1}\right)\right)
\max(A)
\end{align*}
where $A=\left\{\|V_F\| : F \subset G,\,
\mathrm{card}(F)=\mathrm{card}(G)-1\right\}$. Hence:
$$\|V_G\| \leq (2k^2)^{\mathrm{card}(G)}\max(A),$$
since, with $a=(2k^2)^{\mathrm{card}(G)-1} \geq 1$, we have $a+k(k+a)=(1+k)a+k^2 \leq 2k^2a$ because $k \geq 2$.
This ends the proof by induction.
\end{proof}

In the following, we use the multi-index notation. Denote by $\mathbb{Z}_+$ the set of all nonnegative  integers. For $\alpha=(\alpha_1,...,\alpha_n) \in \mathbb{Z}_+^n$, write $|\alpha|=\displaystyle\sum_{i=1}^{n} \alpha_i$. For $j \in \{1,...,n\}$, let $e_j=(\underbrace{0,...,0}_{j-1},1,0,...,0)$.

This last lemma will require the hypothesis $(iv)$ of Theorem \ref{tuples} for the proof to work:

\begin{lemma}\label{Valpha}
Let $X$ be a Banach space.
For all $\alpha \in \mathbb{Z}_+^n$, there exist operators $V_{\alpha} : X \longrightarrow Y$ satisfying:
\begin{enumerate}
    \item[(1)] $V_{(0,...,0)}$ is the operator given in (iv) of Theorem \ref{tuples};
    \item[(2)] $\forall \alpha \in \mathbb{Z}_+^n, \alpha_j=0, T_jV_{\alpha}=0$
    \item[(3)] $\forall \alpha \in \mathbb{Z}_+^n, \alpha_j \geq 1, T_jV_{\alpha}=V_{\alpha-e_j}$
    \item[(4)] $\|V_{\alpha}\| \leq (2k^2)^{n|\alpha|}\|V_{(0,...,0)}\|$.
\end{enumerate}
\end{lemma}

\begin{proof}
Let $V_{(0,...,0)}$ be the operator given in $(iv)$ of Theorem \ref{tuples}. 

We construct the mappings $V_{\alpha}$ by induction on $|\alpha|$. Let $\alpha \in \mathbb{Z}_+^n$, $|\alpha| \geq 1$, and suppose the mappings $V_{\beta} : X \longrightarrow Y$ satisfying (1)-(4) have already been constructed for all $\beta \in \mathbb{Z}_+^n$ satisfying $|\beta|<|\alpha|$. In particular, the mappings $V_{\beta}$ are already constructed for all $\beta \leq \alpha$, $\beta \neq \alpha$.

Let $G=\{j \in \{1,...,n\} : \alpha_j \neq 0\}$ and let $m=\mathrm{card}(G)$. Consider the $m$-tuple of operators $U_j:=T_j\restriction_{\bigcap_{i \notin G} \ker(T_i)} $, $j \in G$. These operators also satisfy the conditions $(i)$, $(ii)$ and $(iii)$ of Theorem \ref{tuples}. Indeed, for $(i)$,  note that $\displaystyle\bigcap_{j \in G} \ker(U_j)$ contains $\displaystyle\bigcap_{j=1}^{n} \ker(T_j)$, for $(ii)$ and $(iii)$ just consider $\widetilde{F}=G^c \cup F \subset \{1,...,n\}$. For $(iii)$, note also that a complement in $Y$ of a subspace of $\displaystyle\bigcap_{i \notin G} \ker(T_i)$ yields a complement in $\displaystyle\bigcap_{i \notin G} \ker(T_i)$.
For $F \subsetneq G$, let 
$$
W_F : X \longrightarrow \displaystyle\bigcap_{i \notin G} \ker(T_i)
$$ 
be defined by $W_F=V_{(\beta_1,...,\beta_n)}$ where $\beta_j=\alpha_j$ if $j \in F$, $\beta_j=\alpha_j-1$ if $j \in G \setminus F$, and $\beta_j=0$ if $j \notin G$. Let us check that the operators $W_F$ satisfy the conditions of Lemma \ref{opVF} with the $U_j$. Recall that $Y_G=\displaystyle\bigcap_{i \notin G} \ker(T_i)$ and for $j \in G, R_j^G:=P_{j,G^c}\widehat{T_j}\restriction_{Y_G}$. Since $\widehat{T_j}y-R_j^Gy \in \ker(T_j)$ (by Lemma \ref{proj}), one has $T_jR_j^Gy=y$. Thus, $R_j^G : Y_G \longrightarrow Y_G$ is a right inverse of $U_j$ with $\|R_j^G\| \leq k$. Let us check that Lemma \ref{opVF} applies to the tuple $(U_j)_{j \in G}$ acting on $Y_G$, with the same constant $k$. 

 (a) Each $U_j$ has the bounded right inverse $R_j^G : Y_G \longrightarrow Y_G$ constructed above, with $\|R_j^G\| \leq k$.

(b) : For $F \subset G$ and $j \in G \setminus F$, the restriction of $P_{j, F \cup G^c}$ to $\displaystyle\bigcap_{i \in F} \ker(U_iU_j)=Y_G \cap \displaystyle\bigcap_{i \in F} \ker(T_iT_j)$ is a projection onto $\displaystyle\bigcap_{i \in F} \ker(U_i)$ whose kernel is contained in $\ker(U_j)$, and its norm is at most $k$; these are the projections of Lemma \ref{proj} for $(U_j)$. 

(c) : By the induction hypothesis, each $W_F, F \subsetneq G$, takes values in $Y_G$ (property (2)) and satisfies $U_jW_F=W_{F \setminus \{j\}}$ for $j \in F$. (property (3)).

Since the proof of Lemma \ref{opVF} only uses (a), (b), (c), and the commutativity of the $U_j$, it yields an operator  $V_{\alpha} : X \longrightarrow \displaystyle\bigcap_{i \notin G} \ker(T_i) \subset Y$ satisfying 
$$T_jV_{\alpha}=W_{G \setminus \{j\}}=V_{\alpha-e_j} \quad (j \in G),$$
$$T_jV_{\alpha}=0 \quad (j \notin G)$$
and
\begin{align*}
\|V_{\alpha}\| &\leq (2k^2)^{\mathrm{card}(G)}\max\{\|W_F\| : F \subset G, \mathrm{card}(F)=\mathrm{card}(G)-1\}\\
&\leq (2k^2)^n\max_{j\in G} \|V_{\alpha -e_j}\| \\
&\leq (2k^2)^n(2k^2)^{n(|\alpha|-1)}\|V_{(0,...,0)}\|  \\
&=(2k^2)^{n |\alpha|} \|V_{(0,...,0)}\|.
\end{align*}
This completes the construction of the operators $V_{\alpha}$.
\end{proof}

One can now prove Theorem \ref{tuples}:

\begin{proof}[Proof of Theorem \ref{tuples}]
    Let $T_1,\ldots,T_n\in \mathcal L(Y)$ be a commuting $n$-tuple
    satisfying the conditions of Theorem \ref{tuples}, and
    let $X$ be the Banach space appearing in hypothesis $(iv)$.
    Let $B_1,\ldots,B_n\in \mathcal L(X)$ be an arbitrary commuting
    $n$-tuple. Choose a scalar $s\in\C\setminus\{0\}$
    sufficiently small and set
    $$
    A_j=sB_j,\qquad
    \gamma=\max\{\|A_j\|:1\leq j\leq n\}.
    $$
    The proof will specify how small $\gamma$ must be.

    Let $V_\alpha:X\to Y$ ($\alpha\in\mathbb Z_+^n$)
    be the operators constructed in Lemma \ref{Valpha},
    and write $V_0=V_{(0,\ldots,0)}$.
    For $\alpha=(\alpha_1,\ldots,\alpha_n)$, set
    $$
    |\alpha|=\alpha_1+\cdots+\alpha_n,
    \qquad
    A^\alpha=A_1^{\alpha_1}\cdots A_n^{\alpha_n}.
    $$
    By the norm estimate in Lemma \ref{Valpha}, we have
    \begin{align*}
        \sum_{\substack{\alpha\in\mathbb Z_+^n\\\alpha\neq0}}
        \|V_\alpha A^\alpha\|
        &\leq
        \|V_0\|
        \sum_{\substack{\alpha\in\mathbb Z_+^n\\\alpha\neq0}}
        (2k^2)^{n|\alpha|}\gamma^{|\alpha|}\\
        &=\|V_0\|\sum_{r=1}^{\infty}
        (2k^2)^{nr}\gamma^r
        \operatorname{card}
        \{\alpha\in\mathbb Z_+^n:|\alpha|=r\}\\
        &=\|V_0\|\sum_{r=1}^{\infty}
        (2k^2)^{nr}\gamma^r\binom{r+n-1}{n-1}\\
        &\leq
        2^{n-1}\|V_0\|\sum_{r=1}^{\infty}
        2^r(2k^2)^{nr}\gamma^r.
    \end{align*}
    This last sum is finite when $2(2k^2)^n\gamma<1$
    and tends to zero as $\gamma\to0$. Choose $s\neq0$
    so small that $2(2k^2)^n\gamma<1$ and $2^{n-1}\|V_0\|\displaystyle\sum_{r=1}^{+\infty} 2^r(2k^2)^{nr}\gamma^r<1$.
    The series $V=\sum_{\alpha\in\mathbb Z_+^n}V_\alpha A^\alpha$ therefore converges absolutely in $\mathcal{L}(X,Y)$ and
    defines a bounded operator satisfying $\|V-V_0\|<1$.
    Since hypothesis $(iv)$ gives $\|V_0x\|\geq\|x\|$
    for every $x\in X$, we obtain
    $$
    \|Vx\|
    \geq\|V_0x\|-\|(V-V_0)x\|
    \geq(1-\|V-V_0\|)\|x\|.
    $$
    Thus $V$ is bounded below. Since $X$ is complete,
    its range $V(X)$ is closed, and $V:X\to V(X)$
    is a bounded isomorphism.

    For every $j\in\{1,\ldots,n\}$ and $x\in X$,
    absolute convergence allows us to apply $T_j$
    termwise and reindex the series. Using the identities
    satisfied by $V_\alpha$ in Lemma \ref{Valpha},
    we obtain
    \begin{align*}
        T_jVx
        &=\sum_{\alpha\in\mathbb Z_+^n}
        T_jV_\alpha A^\alpha x\\
        &=\sum_{\substack{\alpha\in\mathbb Z_+^n\\
                          \alpha_j\geq1}}
        V_{\alpha-e_j}A^\alpha x\\
        &=\sum_{\beta\in\mathbb Z_+^n}
        V_\beta A^{\beta+e_j}x\\
        &=\sum_{\beta\in\mathbb Z_+^n}
        V_\beta A^\beta A_jx
        =VA_jx,
    \end{align*}
    where the last identity uses the commutativity
    of $A_1,\ldots,A_n$.

    Hence $V(X)$ is a closed subspace of $Y$ invariant
    under every $T_j$, and
    $$
    T_j\restriction_{V(X)}\,V=VA_j=V(sB_j),
    \qquad 1\leq j\leq n.
    $$
    Thus $V$ implements a simultaneous similarity between
    $(T_1\restriction_{V(X)},\ldots,T_n\restriction_{V(X)})$
    and $(sB_1,\ldots,sB_n)$, as required.
\end{proof}

The following example is inspired by Example 9 of \cite{Mul}.

\begin{example}
Let $X$ be a Banach space. For $1 \leq p<\infty$, let us define:
$$X_p=\left\{f : \mathbb{Z}_+^n \longrightarrow X : \|f\|^p:=\sum_{\alpha \in \mathbb{Z}_+^n} \|f(\alpha)\|_X^p<\infty\right\}.$$
It is a Banach space for the norm thus defined. Let $(T_1,...,T_n)$ be the tuple of operators in $\mathcal L(X_p)$ defined by:
$$T_j(f)(\alpha)=f(\alpha+e_j) \quad (\alpha \in \mathbb{Z}_+^n)$$
We can exhibit a link between $T_j$ and an operator on the Banach space: 
\begin{align*}Y_p:=\Big\{f \in \mathrm{Hol}(\mathbb{D}^n,X) : f(z_1,...,z_n)=\sum_{(i_1,...,i_n) \in \mathbb{Z}_+^n}z_1^{i_1}...z_n^{i_n}a_{(i_1,...,i_n)}\\ \text{ with } \sum_{(i_1,...,i_n) \in \mathbb{Z}_+^n} \|a_{(i_1,...,i_n)}\|_X^p<\infty \Big\},
\end{align*}
defined by:
$$W_j(f)(z_1,...,z_n)= \frac{f(z_1,...,z_n)-f(z_1,...,z_{j-1},0,z_{j+1},...,z_n)}{z_j} \quad  (z_1,...,z_n) \in \mathbb{D}^n$$
Let us define the operator:
$$\fonction{\Phi}{X_p}{Y_p}{f}{ z \in \mathbb{D}^n \longmapsto\displaystyle\sum_{(i_1,...,i_n) \in \mathbb{Z}_+^n}z_1^{i_1}...z_n^{i_n}f(i_1,...,i_n)}$$
One can easily check that $\Phi$ is an isometric isomorphism and that $W_j=\Phi T_j\Phi^{-1}$.

Let us check that the operators $(T_1,...,T_n)$ satisfy the hypotheses of Theorem \ref{tuples}.
\begin{enumerate}
    \item[$(i)$] $\displaystyle\bigcap_{j=1}^{n} \ker(T_j)=\displaystyle\bigcap_{j=1}^{n} \ker(W_j)$ is the set of all functions in $Y_p$ which are constant with respect to all of its variables. So $\displaystyle\bigcap_{j=1}^{n} \ker(T_j)$ is isomorphic to $X$ and so is infinite-dimensional.

    \item[$(ii)$] Let $F \subset \{1,...,n\}$ and $j \in \{1,...,n\} \setminus F$. Let $f \in \displaystyle\bigcap_{i \in F} \ker(T_i)$. Let $g \in X_p$ defined by:
    $$g(\alpha_1,...,\alpha_n)=\begin{cases}f(\alpha_1,...,\alpha_{j-1},\alpha_j-1,\alpha_{j+1},...,\alpha_n) \text{ if } \alpha_j\geq 1\\
    0 \text{ if } \alpha_j=0.
    \end{cases}$$
    Then, $g \in \displaystyle\bigcap_{i \in F} \ker(T_i)$ and $T_j(g)=f$.

    \item[$(iii)$] Let $F \subset \{1,...,n\}$ and $j \in \{1,...,n\} \setminus F$. Let $A=\displaystyle\bigcap_{i \in F} \ker(T_i) \cap \ker(T_j)$. With
    $$E=\{\alpha \in \mathbb{Z}_+^n : \alpha_k=0 \text{ for all } k \in F \cup \{j\}\},$$
    one has $A=\{f \in X_p : f=0 \text{ outside } E\}$ and $Pf=\mathbb{1}_Ef$ is a projection of norm 1 onto $A$. Hence, $A$ is complemented.
    
    \item[$(iv)$] We saw in (1) that $\displaystyle\bigcap_{j=1}^{n} \ker(T_j)$ is isomorphic to $X$. So there exists an operator $V_{(0,...,0)} : X \longrightarrow \displaystyle\bigcap_{j=1}^{n} \ker(T_j)$ which is bounded below, hence, we get the existence of an operator $V_{(0,...,0)} : X \longrightarrow \displaystyle\bigcap_{j=1}^{n} \ker(T_j)$ which satisfies $(iv)$.
\end{enumerate}
Thus, thanks to Theorem \ref{tuples}, the $n$-tuple $(T_1,...,T_n)$ is universal for all commuting tuples on $X$. In particular, for $X=\ell_p$, one has $X_p \simeq \ell_p$, which gives a commuting $n$-tuple on $\ell_p$ universal for commuting $n$-tuples of operators on $\ell_p$.
\end{example}

For $n=2$, the hypothesis (ii) of Theorem \ref{tuples} can be simplified:

\begin{lemma}
Let $T_1,T_2 \in \mathcal L(Y)$ be commuting surjective operators. The following statements are equivalent:
\begin{enumerate}
    \item[(1)] $\ker(T_1T_2)=\ker(T_1) +\ker(T_2)$;
    \item[(2)] $T_1(\ker(T_2))=\ker(T_2)$;
    \item[(3)] $T_2(\ker(T_1))=\ker(T_1)$. 
\end{enumerate}
\end{lemma}

\begin{proof}
(1) $\implies$ (2): Clearly, $T_1(\ker(T_2)) \subset \ker(T_2)$ since $T_1$ and $T_2$ are commuting. Let $x \in \ker(T_2)$. Since $T_1$ is surjective, there exists $y \in Y$ such that $T_1(y)=x$. Thus $T_1T_2(y)=0$ and by the assumption, $y=y_1+y_2$ for some $y_1 \in \ker(T_1)$ and $y_2 \in \ker(T_2)$. Then, 
$$T_1(y-y_1)=T_1(y)=x \text{ and } T_2(y-y_1)=T_2(y_2)=0$$
So $x \in T_1(\ker(T_2))$.

(2) $\implies$ (1): The inclusion $\ker(T_1) + \ker(T_2) \subset\ker(T_1T_2)$ is always true.

Let $x \in \ker(T_1T_2)$. Then $T_1(x) \in \ker(T_2)$, and so, by assumption (2), there exists $y \in \ker(T_2)$ such that $T_1(y)=T_1(x)$. Hence, 
$x-y \in \ker(T_1)$, and $x=x-y+y \in \ker(T_1)+\ker(T_2)$.

The equivalence (1) $\iff$ (3) follows by symmetry.
\end{proof}

This yields an immediate corollary for $n=2$:

\begin{corollary}
Let $X,Y$ be Banach spaces. Let $T_1,T_2 \in \mathcal L(Y)$ be commuting surjective operators satisfying 
\begin{enumerate}
    \item[(1)] $\dim(\ker(T_1) \cap \ker(T_2))=\infty$;
    \item[(2)] $\ker(T_1T_2)=\ker(T_1)+\ker(T_2)$;
    \item[(3)] $\ker(T_1) \cap \ker(T_2)$ is complemented in $Y$;
    \item[(4)]
    $\ker(T_1)$ and $\ker(T_2)$ are complemented in $Y$ (equivalently, $T_1$ and $T_2$ have bounded right inverses);
    \item[(5)] There exists an operator $V : X \longrightarrow \ker(T_1) \cap \ker(T_2)$ such that for all $x \in X, \|V(x)\| \geq \|x\|$. 
\end{enumerate}
Then, the pair $(T_1,T_2)$ is universal for all commuting pairs on $X$.
\end{corollary}

\subsection{Non-commuting \texorpdfstring{$n$}{n}-tuples}

In his article \cite{Mul}, V. Müller also provided a theorem which gives sufficient conditions for a non-commuting $n$-tuple of operators to be universal. Let us generalize his result to a Banach space setting.

\begin{definition}\label{defntuples1}
We say that an $n$-tuple $T_1,...,T_n \in \mathcal L(Y)$ of operators is universal for all $n$-tuples on $X$ if it has the following property: for each $n$-tuple $A_1,...,A_n \in \mathcal L(X)$ there exist a constant $c \neq 0$ and a closed subspace $M \subset Y$ invariant for $T_1,...,T_n$ such that the $n$-tuples $(cA_1,...,cA_n)$ and $(T_1\restriction_M,...,T_n\restriction_M)$ are similar.
\end{definition}

\begin{theorem}\label{tuplesnc}
Let $X,Y$ be complex Banach spaces and let $T_1,...,T_n \in \mathcal L(Y)$ satisfy the following properties:
\begin{enumerate}
    \item[(1)] $\dim\left(\displaystyle\bigcap_{j=1}^{n} \ker(T_j)\right)=\infty$;
    \item[(2)] $T_j\left(\displaystyle\bigcap_{i \neq j} \ker(T_i)\right)=Y$ for all $j \in \{1,...,n\}$;
    \item[(3)] There exists an operator denoted $V_{\theta} : X \longrightarrow \displaystyle\bigcap_{j=1}^{n} \ker(T_j)$ such that $\|V_{\theta}(x)\| \geq \|x\|$ for all $x \in X$. 
    \item[(4)] $\displaystyle \bigcap_{i=1}^n \ker(T_i)$ is complemented in $Y$.
\end{enumerate}
Then the $n$-tuple $(T_1,...,T_n)$ is universal for all $n$-tuples on $X$.
\end{theorem}

\begin{remark}
Again, $(3) \implies (1)$, and the conditions $(1)$ and $(3)$ are necessary. 

We can weaken the assumption $(4)$, by replacing it by 
$$(4') ~ \bigcap_{i=1}^{n} \ker(T_i) \text{ is complemented in } \bigcap_{i \neq j} \ker(T_i) \text{ for every } j,$$
at the cost of a less transparent statement.
\end{remark}

\begin{proof}
For $j \in \{1,...,n\}$, let 
$\widehat{T_j} : Y \longrightarrow \displaystyle\bigcap_{i \neq j} \ker(T_i)$ be a bounded right inverse of the restriction of $T_j$ to the subspace $\displaystyle\bigcap_{i \neq j} \ker(T_i)$. The existence of $\widehat{T_j}$ is ensured by assumptions $(2)$ and $(4)$. Let $k=\max\{\|\widehat{T_j}\| : j \in \{1,...,n\}\}$. 

For $r \geq 0$, let $F_r=\{(\alpha_r,\alpha_{r-1},...,\alpha_1) : \alpha_j \in \{1,...,n\}\}$. Clearly, $\mathrm{card}(F_r)=n^r$. Let $\mathcal{F}=\displaystyle\bigcup_{r=0}^{+\infty} F_r$. For $r=0$, the only element of $F_0$ will be denoted $\theta$.

Let $A_1,...,A_n \in \mathcal L(X)$. Let $V_{\theta} : X \longrightarrow \displaystyle\bigcap_{j=1}^{n} \ker(T_j)$ be the operator defined in hypothesis (3). Choose $s \neq 0$ such that $\gamma:=|s|\underset{1 \leq i \leq n}{\max}\|A_i\|$ satisfies $\gamma nk<1$ and $\frac{\|V_{\theta}\|\gamma nk}{1-\gamma nk}<1$.  

For $(\alpha_r,...,\alpha_1) \in \mathcal{F}$, define the operator $V_{(\alpha_r,...,\alpha_1)} : X \longrightarrow Y$ by:
$$V_{(\alpha_r,...,\alpha_1)}=\widehat{T_{\alpha_r}}...\widehat{T_{\alpha_1}}V_{\theta}.$$
Then, $T_jV_{(\alpha_r,...,\alpha_1)}=0$ if $j \neq \alpha_r$ since $\widehat{T_{\alpha_r}}$ takes its values in particular in $\ker(T_j)$. If $\alpha_r=j$, then $T_jV_{(\alpha_r,...,\alpha_1)}=V_{(\alpha_{r-1},...,\alpha_1)}$. Moreover, $\|V_{(\alpha_r,...,\alpha_1)}\| \leq k^r\|V_{\theta}\|$. 
Define $V : X \longrightarrow Y$ by:
$$V(x)=\sum_{(\alpha_r,...,\alpha_1) \in \mathcal{F}} V_{(\alpha_r,...,\alpha_1)}s^rA_{\alpha_1}...A_{\alpha_r}(x) \quad (x \in X)$$
We have
$$\|V-V_{\theta}\| \leq \sum_{(\alpha_r,...,\alpha_1) \in \mathcal{F} \setminus F_0} \|V_{(\alpha_r,...,\alpha_1)}\|\gamma^r \leq \|V_{\theta}\|\sum_{r=1}^{+\infty} k^rn^r\gamma^r=\|V_{\theta}\|\frac{\gamma nk}{1-\gamma nk}<1$$
So $V$ is a bounded linear operator. Since $V_{\theta}$ is bounded below by 1, $V$ is bounded below and its range $M:=V(X)$ is a closed subspace of $Y$.

For all $j \in \{1,...,n\}$ and $x \in X$ we have 
$$VA_j(x)=\sum_{(\alpha_r,...,\alpha_1) \in \mathcal{F}} V_{(\alpha_r,...,\alpha_1)}s^rA_{\alpha_1}...A_{\alpha_r}A_j(x)$$
and
\begin{align*}T_jV(x)&=\sum_{(\alpha_r,...,\alpha_1) \in \mathcal{F}} T_jV_{(\alpha_r,...,\alpha_1)}s^rA_{\alpha_1}...A_{\alpha_r}(x)\\&=\sum_{(\alpha_{r-1},...,\alpha_1)} V_{(\alpha_{r-1},...,\alpha_1)}s^rA_{\alpha_1}...A_{\alpha_{r-1}}A_j(x)
\end{align*}
So $T_jV=sVA_j$. Hence $M$ is a closed subspace of $Y$ invariant for all $T_1,...,T_n$, and the $n$-tuples $(sA_1,...,sA_n)$ and $(T_1\restriction_M,...,T_n\restriction_M)$ are similar.
\end{proof}

\begin{example}
In light of these theorems, one can provide an example of a universal pair of operators in $\ell_p$, $1 \leq p<\infty$. Let $(e_n)$ be the canonical basis of $\ell_p$. Define the following operators:
$$\begin{array}{lllllcc}
    T_1 : & \ell_p & \longrightarrow & \ell_p \\
     & e_{4n} & \longmapsto & 0\\
     & e_{4n+1} & \longmapsto & e_n\\
     & e_{4n+2} & \longmapsto & 0\\
     & e_{4n+3} & \longmapsto & 0
\end{array} \quad \text{ and } \begin{array}{lllllcc}
    T_2 : & \ell_p & \longrightarrow & \ell_p \\
     & e_{4n} & \longmapsto & 0\\
     & e_{4n+1} & \longmapsto & 0\\
     & e_{4n+2} & \longmapsto & e_n\\
     & e_{4n+3} & \longmapsto & 0
\end{array}$$
Then, $\ker(T_1) \cap \ker(T_2)=\overline{\mathrm{span}}\{e_{4n},e_{4n+3}\}$ is infinite-dimensional and it is complemented by the closed linear span of the remaining coordinates. We also have that for $j=1,2$ and $i \neq j$, $T_j\left(\ker(T_i)\right)=\ell_p$. Moreover, the operator:
$$\begin{array}{lllllcc}
    V : & \ell_p & \longrightarrow & \overline{\mathrm{span}\{e_{4n},e_{4n+3}\}} \\
     & e_{4n} & \longmapsto & e_{8n}\\
     & e_{4n+1} & \longmapsto & e_{8n+4}\\
     & e_{4n+2} & \longmapsto & e_{8n+3}\\
     & e_{4n+3} & \longmapsto & e_{8n+7}
\end{array}$$
is an isometry. 
Thus, the pair $(T_1,T_2)$ satisfies the hypothesis of Theorem \ref{tuplesnc} (with $X=Y=\ell_p$), so it is universal for all pairs of operators on $\ell_p$.
\end{example}

\subsection{About the essential spectrum of universal \texorpdfstring{$n$}{n}-tuples of operators}

This last part shows that the essential spectrum of a universal $n$-tuple of operators, commuting or not, is large. For a $n$-tuple $T=(T_1,...,T_n) \in \mathcal{L}(Y)^n$ and $a=(a_1,...,a_n) \in \mathbb{C}^n$, we write 
$$\delta_{T-a} : Y \longrightarrow Y^n, \quad y \longmapsto ((T_1-a_1)y,...,(T_n-a_n)y).$$

\begin{definition}\label{defntuples2}
Let $Y$ be a Banach space and $T=(T_1,...,T_n) \in \mathcal{L}(Y)^n$. The \emph{joint approximate point spectrum} and the \emph{joint essential approximate point spectrum} of $T$ are 
$$\sigma_{\pi}(T)=\left\{a \in \mathbb{C}^n : \exists (y_k)_{k \geq 1} \subset Y, \|y_k\|=1, T_iy_k-a_iy_k \underset{k \rightarrow +\infty}{\longrightarrow} 0 ~ \forall i \in \{1,...,n\}\right\}$$
\begin{align*}
\sigma_{\pi e}(T)=\Big\{a \in \mathbb{C}^n : \exists (y_k)_{k \geq 1} \subset Y,  \|y_k\|=1, (y_k) \text{ has no convergent subsequence, }\\
T_ix_k-a_ix_k \underset{k \rightarrow +\infty}{\longrightarrow} 0, \forall i \in \{1,...,n\} \Big\}.
\end{align*}
\end{definition}

\begin{remark}\label{rem:spec}
    \begin{enumerate}
        \item A bounded operator $S$ between Banach spaces has finite-dimensional kernel and closed range if and only if every bounded sequence $(y_k)$ with $Sy_k \underset{k \rightarrow +\infty}{\longrightarrow} 0$ has a convergent subsequence. Hence $a \notin \sigma_{\pi e}(T)$ if and only if $\delta_{T-a}$ is upper semi-Fredholm (see \cite{Muller}). For $n=1$, $\sigma_{\pi e}(S)=\sigma_e^+(S,Y)$ is the set introduced in Section \ref{spectrum}.  

        \item The condition "no convergent subsequence" cannot be dropped. For the backward shift $B$ on $\ell_p$, every $a \in \mathbb{D}$ is an eigenvalue of multiplicity one and $B-a I=B(I-aS)$ is surjective ($S$ is the forward shift), so $B-aI$ is Fredholm. Therefore, $\sigma_{\pi}(B)=\overline{\mathbb{D}}$, whereas $\sigma_{\pi e}(B)=\mathbb{T}$.

        \item If $T$ is commuting, $\sigma_{\pi e}(T)$ is contained in the Taylor essential spectrum $\sigma_{Te}(T)$, i.e. the set of $a$ for which the Koszul complex of $T-a$ is not Fredholm, and in the essential (Calkin) spectrum of $T$ in the sense of Harte \cite{Curto}, \cite {Muller}. Indeed, the Koszul complex of $T-a$ is Fredholm, its first map, which is $\delta_{T-a}$, has closed range and finite-dimensional kernel.

        \item If $\displaystyle\bigcap_{i=1}^{n} \ker(T_i-a_iI)$ is infinite-dimensional, then $a \in \sigma_{\pi e}(T)$. Recall \emph{Riesz's lemma}: if $F$ is a closed subspace of a normed space $E$ with $F \neq E$ and $0<\theta<1$, there exists $x \in E$ with $\|x\|=1$ and $\mathrm{dist}(x,F) \geq \theta$. Applying it with $\theta=\frac{1}{2}$ to the finite-dimensional, hence closed, subspaces $\mathrm{span}(y_1,...,y_k)$ of the infinite-dimensional closed subspace $E=\displaystyle\bigcap_{i=1}^{n} \ker(T_i-a_iI)$, we obtain inductively unit vectors $y_k \in E$ with $\|y_k-y_{\ell}\| \geq \frac{1}{2}$ for $k \neq \ell$. They have no convergent subsequence and $T_iy_k=a_iy_k$ for all $i$ and $k$.
    \end{enumerate}
\end{remark}

\begin{theorem}
    Let $1 \leq p<\infty$ and let $T=(T_1,...,T_n)\in \mathcal{L}(\ell_p)^n$ be universal for commuting $n$-tuples on $\ell_p$ (Definition \ref{univtuples}) or universal for all $n$-tuples on $\ell_p$ (Definition \ref{defntuples1}). Let $\alpha_1,...,\alpha_{n-1} \in \mathbb{C}$ and set $\alpha_0=1$. There exists $r>0$ such that for every $b \in r\mathbb{D}$ the joint eigenspace 
    $$\bigcap_{i=1}^{n} \ker(T_i-\alpha_{i-1}bI)$$
    is infinite-dimensional. In particular, 
    $$r(1,\alpha_1,...,\alpha_{n-1})\mathbb{D} \subset \sigma_{\pi e}(T) \subset \sigma_{\pi}(T),$$
    and, if $T$ is commuting, $r(1,\alpha_1,...,\alpha_{n-1})\mathbb{D} \subset \sigma_{Te}(T)$.
\end{theorem}

\begin{proof}
Consider again $X_0=\displaystyle{\left(\bigoplus_{\mathbb{N} \cup \{0\}} \ell_p\right)_{\ell_p}}$ and recall that $X_0 \simeq \ell_p$.

Let $B_{\infty}(x_0,x_1,x_2,...)=(x_1,x_2,...)$ be the backward shift on $X_0$. For $a \in \mathbb{D}$ and $u \in \ell_p$, put $e_{a,u}=(u,au,a^2u,...) \in X_0$. Then $B_{\infty}e_{a,u}=ae_{a,u}$ and $\|e_{a,u}\|=(1-|a|^p)^{-\frac{1}{p}}\|u\|$, so $u \longmapsto e_{a,u}$ maps $\ell_p$ isomorphically onto a subspace of $\ker(B_{\infty}-aI)$. In particular, $\ker(B_{\infty}-aI)$ is infinite-dimensional for every $a \in \mathbb{D}$. 

The $n$-tuple $A=(B_{\infty},\alpha_1 B_{\infty},...,\alpha_{n-1}B_{\infty})$ on $X_0 \simeq \ell_p$ is commuting. By universality, there exist $c \neq 0$, a closed subspace $M \subset \ell_p$ invariant under $T_1,...,T_n$, and an isomorphism $J : X_0 \longrightarrow M$ such that $T_iJ=cJ\alpha_{i-1}B_{\infty}$ for $i \in \{1,...,n\}$. The constant $c$ depends on $\alpha_1,...,\alpha_{n-1}$. For $a \in \mathbb{D}$ and $y \in \ker(B_{\infty}-aI)$, we get
$$T_iJy=c\alpha_{i-1}JB_{\infty}y=\alpha_{i-1}(ca)Jy \quad (i \in \{1,...,n\}),$$
that is, $J(\ker(B_{\infty}-aI)) \subset \displaystyle\bigcap_{i=1}^{n} \ker(T_i-\alpha_{i-1}(ca)I)$. Since $J$ is injective, this joint eigenspace is infinite-dimensional. As $a$ runs over $\mathbb{D}$, $b=ca$ runs over $c\mathbb{D}=|c|\mathbb{D}$, so $r=|c|$ satisfies the first assertion. The inclusion follows from Remark \ref{rem:spec} (4) and (3).
\end{proof}

\begin{remark}
    The proof gives more than membership in an essential spectrum: every point of $r(1,\alpha_1,...,\alpha_{n-1})\mathbb{D}$ is a joint eigenvalue of $T$ of infinite multiplicity. If one prefers an explicit sequence, $y_k=\frac{J(e_{a,e_k})}{\|J(e_{a,e_k})\|}$, where $(e_k)$ is the unit vector basis of $\ell_p$, has no convergent subsequence. Indeed, if $(y_{k_j})$ converged, a further subsequence of $(J(e_{a,e_{k_j}}))$ would converge since $\|J^{-1}\|^{-1}\|e_{a,e_k}\| \leq \|J(e_{a,e_k})\| \leq \|J\|\|e_{a,e_k}\|$. Then, $(e_{a,e_{k_j}})$ would converge, which is impossible because $\|e_{a,e_k}-e_{a,e_{\ell}}\|=2^{\frac{1}{p}}(1-|a|^p)^{-\frac{1}{p}}$ for $k \neq \ell$.
\end{remark}

\section{Acknowledgments}

The author wants to express sincere thanks to Isabelle Chalendar and Colin Petitjean for many useful discussions about the subject of this work, and for their careful reading of the paper.

\section{AI Statement}

The author declares that he did not use any AI to do the research work nor to write the article.



\begin{thebibliography}{99}

\bibitem{AK} F. Albiac, N.J. Kalton \emph{Topics in Banach Space Theory}, Graduate Texts in Mathematics 233. Berlin: Springer (ISBN 0-387-28141-X/hbk). xi, 373 p. (2006).

\bibitem{AndroulakisSchlumprecht2003}
G. Androulakis, T. Schlumprecht,
\emph{The Banach space $S$ is complementably minimal and subsequentially prime},
Studia Math. 156, No. 3, 227--242 (2003).

\bibitem{AH2011}
S.A. Argyros and R.G. Haydon,
\emph{A hereditarily indecomposable ${\mathcal L}_{\infty}$-space that solves the scalar-plus-compact problem},
Acta Math., 206(1):1--54, 2011.

\bibitem{Bourgain} J. Bourgain, \emph{On separable Banach spaces, universal for all separable reflexive spaces}, Proc. Amer. Math. Soc. \textbf{79} (1980), 241--246.

\bibitem{Caradus1969}
S.R. Caradus,
\emph{Universal operators and invariant subspaces},
Proc. Amer. Math. Soc., 23:526--527, 1969.

\bibitem{CasazzaJohnsonTzafriri1984}
P. G. Casazza, W. B. Johnson, L. Tzafriri,
\emph{On Tsirelson's space},
Israel J. Math. 47, No. 2--3, 81--98 (1984).

\bibitem{CP2011}
I.~Chalendar and J.~R.~Partington,
\emph{Modern Approaches to the Invariant-Subspace Problem}, volume 188 of Cambridge Tracts in Mathematics.
Cambridge University Press, Cambridge, 2011.

\bibitem{CGGintro}
C.C. Cowen and E.A. Gallardo-Guti\'errez.
\emph{An introduction to Rota's universal operators: properties, old and new examples and future issues},
Concr. Oper., 3:43--51, 2016.

\bibitem{CGGnewproof}
C.C. Cowen and E.A. Gallardo-Guti\'errez.
\emph{A new proof of a Nordgren, Rosenthal and Wintrobe theorem on universal operators.}
In Problems and Recent Methods in Operator Theory (Memphis, TN, 2015),
97--102, Contemp. Math. 687, Amer. Math. Soc., 2017.

\bibitem{Curto} R. E. Curto, \emph{Fredholm and invertible $n$-tuples of operators. The deformation problem}, Trans. Amer. Math. Soc. \textbf{266} (1981), 129--159.

\bibitem{Enflo}
P. Enflo,
\emph{On the invariant subspace problem in Banach spaces},
S\'eminaire Analyse fonctionnelle (Maurey--Schwartz), 1975--1976, exp.\ no.\ A11, 1--6.

\bibitem{Gowers} W.T. Gowers,
\emph{A new dichotomy for Banach spaces},
Geom. Funct. Anal. 6 (1996), 1083-1093.

\bibitem{GM}
W.T. Gowers, B. Maurey,
\emph{The unconditional basic sequence problem},
J. Am. Math. Soc. 6, No. 4, 851-874 (1993).

\bibitem{KT} A. Ryszard Komorowski, N. Tomczak-Jaegermann, \emph{Banach spaces without local unconditional structure}, Isr. J. Math. 89, No. 1-3, 205-226 (1995).

\bibitem{Lindenstrauss1967}
J. Lindenstrauss,
\emph{On complemented subspaces of $m$}
Israel J. Math. 5, 153--156 (1967).

\bibitem{LT1977}
J. Lindenstrauss and L. Tzafriri. \emph{Classical Banach Spaces. I},
Springer-Verlag, Berlin--New York, 1977.

\bibitem{muj} 
J. Mujica, \emph{Complex analysis in Banach spaces. Holomorphic functions and domains of holomorphy in finite and infinite dimensions}, North-Holland Mathematics Studies, 120. Notas de Matemática, 107. Amsterdam/New York/Oxford: North-Holland. XI, 434 p. (1986). 

\bibitem{Muller} V. Müller, \emph{Spectral Theory of Linear Operators and Spectral Systems in Banach Algebras, 2nd ed.}, Operator Theory: Advances and Applications 139. Basel: Birkhäuser (ISBN 978-3-7643-8264-3/hbk). viii, 439 p. (2007). 

\bibitem{Mul}
V. Müller

\emph{Universal $n$-tuples of operators}, Mathematical Proceedings of the Royal Irish Academy 113A (2013), 143-150; doi:10.3318/PRIA.2013.113.13

\bibitem{NRW1987}
E. Nordgren, P. Rosenthal, and F.S. Wintrobe.
\emph{Invertible composition operators on $H^p$},
J. Funct. Anal., 73(2):324--344, 1987.

\bibitem{Read}
C.J. Read.
\emph{A solution to the invariant subspace problem},
Bull. London Math. Soc., 16(4):337--401, 1984.

\bibitem{Read2} C.J. Read,
\emph{A solution to the invariant subspace problem on $\ell_1$},
Bull. London Math Soc 17 (1985) n°4, 305-317.

\bibitem{Rota1959}
G.C. Rota.
\emph{Note on the invariant subspaces of linear operators},
Rend. Circ. Mat. Palermo (2), 8:182--184, 1959.

\bibitem{Rota1960}
G.C. Rota,
\emph{On models for linear operators},
Comm. Pure Appl. Math. 13 (1960), 469--472

\bibitem{schro}
R. Schroderus, H-O Tylli,
\emph{On universal operators and universal pairs},
Proc. Edinb. Math. Soc., II. Ser. 61, No. 3, 891-908 (2018).

\bibitem{Szlenk} W. Szlenk, \emph{The non-existence of a separable reflexive Banach space universal for all separable reflexive Banach spaces}, Studia Math. \textbf{30} (1968), 53--61.

\end{thebibliography}
\end{document}